\documentclass[12pt, reqno]{amsart}
\usepackage{amsmath, amsthm, amscd, amsfonts, amssymb, graphicx, color}
\usepackage[bookmarksnumbered, colorlinks, plainpages]{hyperref}
\hypersetup{colorlinks=true,linkcolor=red, anchorcolor=green, citecolor=cyan, urlcolor=red, filecolor=magenta, pdftoolbar=true}

\newtheorem{theorem}{Theorem}[section]
\newtheorem{lemma}[theorem]{Lemma}

\newtheorem{corollary}[theorem]{Corollary}
\theoremstyle{definition}

\theoremstyle{remark}

\numberwithin{equation}{section}

\begin{document}
\setcounter{page}{1}

\title[Riemann surfaces and $AF$-algebras]{Riemann surfaces and $AF$-algebras}

\author[I.Nikolaev]
{Igor Nikolaev$^1$}

\address{$^{1}$ The Fields Institute for Research in Mathematical Sciences, Toronto, ON, Canada.}
\email{\textcolor[rgb]{0.00,0.00,0.84}{igor.v.nikolaev@gmail.com}}

\dedicatory{In memory of Ola Bratteli}



\date{Received:  August 14, 2015; Revised: yyyyyy; Accepted: zzzzzz.}

\begin{abstract}
For a generic  set in the Teichm\"uller space, we construct a 
covariant functor with the range in a category of the
$AF$-algebras;  the functor maps isomorphic Riemann surfaces to
the stably isomorphic $AF$-algebras.   In the special case of genus one,   
one gets  a  functor   between the category of  complex tori and the
 Effros-Shen  algebras. 

\end{abstract} \maketitle

\section{Introduction}

The aim of our  note  is  a functor from the set of 
generic Riemann surfaces to a category of  the  operator algebras,
known as the $AF$-algebras;  for the sake of clarity,   consider  the simplest
 example.
 We shall write $\{\Lambda_{\tau}={\Bbb Z}+{\Bbb Z}\tau ~|~ \Im ~(\tau)>0\}$
 to denote a lattice in the complex plane ${\Bbb C}$.  
 Let ${\Bbb C}/\Lambda_{\tau}$ be a  complex torus corresponding to $\Lambda_{\tau}$,
 i.e. the  Riemann surface of genus $g=1$.   
 We shall write $\{ {\Bbb A}_{\theta} ~|~ \theta\in {\Bbb R}\}$ to denote an $AF$-algebra
 defined  by the inductive  limit of positive isomorphisms:  
\begin{equation}
{\Bbb Z}^2
\buildrel\rm
\left(
\begin{matrix}
a_0 & 1\cr 1 & 0
\end{matrix}
\right)
\over\longrightarrow 
{\Bbb Z}^2
\buildrel\rm
\left(
\begin{matrix}
a_1 & 1\cr  1 & 0
\end{matrix}
\right)
\over\longrightarrow
{\Bbb Z}^2
\buildrel\rm
\left(
\begin{matrix}
a_2 & 1\cr  1 & 0
\end{matrix}
\right)
\over\longrightarrow
\dots,
\end{equation}
where the regular continued fraction $[a_0, a_1, a_2, \dots]$ 
converges to $\theta$;  we refer the reader to  [Bratteli 1972]   \cite{Bra1}
for a  definition of the $AF$-algebras and [Effros \& Shen 1980]   \cite{EfSh2}
for the properties of  algebra ${\Bbb A}_{\theta}$ (the {\it Effros-Shen algebra}).    
Recall  that  ${\Bbb C}/\Lambda_{\tau}$ and  
${\Bbb C}/\Lambda_{\tau'}$ are isomorphic complex tori, if and only if,  $\tau'={a\tau +b\over c\tau +d}$
for some integers $a,b,c$ and $d$,   such that $ad-bc=\pm 1$;  here an isomorphism means a conformal
map between the Riemann surfaces ${\Bbb C}/\Lambda_{\tau}$ and  ${\Bbb C}/\Lambda_{\tau'}$.  
It is a deep and  amazing fact,  that the same is true of the Effros-Shen algebras.  Namely, 
recall that the $C^*$-algebras ${\Bbb A}$ and ${\Bbb A}'$ are called {\it stably isomorphic}
(Morita equivalent), if ${\Bbb A}\otimes {\Bbb K}\cong {\Bbb A}'\otimes {\Bbb K}$,   where ${\Bbb  K}$
is the $C^*$-algebra of compact operators on a Hilbert space $H$.  
It is known,   that  the Effros-Shen algebras ${\Bbb A}_{\theta}$ and ${\Bbb A}_{\theta'}$ are stably isomorphic,
if and only if,    $\theta'={a\theta +b\over c\theta +d}$   for some integers $a,b,c$ and $d$, such that $ad-bc=\pm 1$,
see  e.g.  [Effros \& Shen 1980]   \cite{EfSh2},  pp. 199-201.   
(A  relation  between complex tori and   continued fractions  was already known to
 [Klein 1896]  \cite{Kle1}.)   
 One may wonder,   if there exists a functor from the category of complex tori  (Riemann surfaces,  resp.)  
to the category of Effros-Shen algebras ($AF$-algebras, resp.),  such that isomorphisms between the Riemann surfaces 
generate stable isomorphisms between the corresponding $AF$-algebras.

In the present paper we construct a covariant
functor $F$ from a generic set of the Riemann surfaces of genus $g\ge 1$
to a category of the so-called toric $AF$-algebras (to be specified below);
the functor maps isomorphic Riemann surfaces to the stably isomorphic 
toric $AF$-algebras (Theorem \ref{thm1}).    
To formulate our results, 
denote by $T(g)$ the Teichm\"uller space of genus $g\ge 1$ and let $S\in T(g)$
be a Riemann surface.
Let  $q\in H^0(S, \Omega^{\otimes 2})$ be a holomorphic quadratic 
differential on the Riemann surface $S$, such that all zeroes
of $q$ are simple, see  [Strebel 1984]  \cite{S}. 
By $\widetilde S$ we understand a double 
cover of $S$ ramified over the zeroes of $q$. Note that there
is an involution on the homology groups $H_*(\widetilde S)$
induced by the covering map $\widetilde S\to S$. 
Let $H_1^{odd}(\widetilde S)$ be the odd part of the first (integral)
homology of $\widetilde S$ with respect to this involution, relatively the  zeroes of $q$. 
By the formulas for the relative homology one gets 
$H_1^{odd}(\widetilde S)\cong {\Bbb Z}^n$, where 
$n=6g-6$ if $g\ge 2$ and $n=2$ if $g=1$. 
It is known, that 
\begin{equation}\label{eq-2}
Hom~(H_1^{odd}(\widetilde S); {\Bbb R})-\{0\} \cong T(g),
\end{equation}
where $0$ is the zero homomorphism  [Hubbard \& Masur  1979]   \cite{HuMa1}. 
Fix a basis in  homology group $H_1^{odd}(\widetilde S)$ and, in view of (\ref{eq-2}),  
 denote  by $(\lambda_1,\dots,\lambda_n)$ its image in ${\Bbb R}$ such that $\lambda_1\ne 0$;
    let  $\theta=(\theta_1,\dots,\theta_{n-1})$ be 
a vector with the coordinates $\theta_i=\lambda_{i-1}/\lambda_1$.
We shall consider the following Jacobi-Perron continued fraction:
\begin{equation}\label{eq-3}
\left(
\begin{matrix}
1\cr \theta
\end{matrix}
\right)=
\lim_{k\to\infty} 
\left(
\begin{matrix}
0 & 1\cr I & b_1
\end{matrix}
\right)\dots
\left(
\begin{matrix}
0 & 1\cr I & b_k
\end{matrix}
\right)
\left(
\begin{matrix}
0\cr {\Bbb I}
\end{matrix}
\right),
\end{equation}
where $b_i=(b^{(i)}_1,\dots, b^{(i)}_{n})^T$ is a vector of the  non-negative integers,  
$I$ the unit matrix and ${\Bbb I}=(0,\dots, 0, 1)^T$;   we refer the reader to
[Bernstein 1971]  \cite{B} for the theory of such fractions.  
Finally, consider an $AF$-algebra   defined by the following inductive limit of
 positive isomorphisms:
\begin{equation}
{\Bbb Z}^n
\buildrel\rm
\left(
\begin{matrix}
0 & 1\cr I & b_1
\end{matrix}
\right)
\over\longrightarrow 
{\Bbb Z}^n
\buildrel\rm
\left(
\begin{matrix}
0 & 1\cr I & b_2
\end{matrix}
\right)
\over\longrightarrow
{\Bbb Z}^n
\buildrel\rm
\left(
\begin{matrix}
0 & 1\cr I & b_3
\end{matrix}
\right)
\over\longrightarrow
\dots,
\end{equation}
see e.g. [Effros 1981]   \cite{E}. 
We shall denote such an algebra by  ${\Bbb A}_{\theta}$
and refer to ${\Bbb A}_{\theta}$ as a {\it toric $AF$-algebra}.  
Notice that if $g=1$,  then the Jacobi-Perron continued fraction
coincides with a regular continued fraction;  thus for $g=1$   
the toric $AF$-algebra is isomorphic  to  an  Effros-Shen
algebra  ${\Bbb A}_{\theta}$,  hence our notation.

Let $F: T(g)\to \{$toric  $AF$-algebras$\}$  be a map acting by
the formula $(\lambda_1,\dots, \lambda_n)\mapsto {\Bbb A}_{\theta}$,
where $\theta=(\theta_1,\dots, \theta_n)$.   Let $V$ be the maximal 
subset of $T(g)$,  such that every  Riemann surface $S\in V$ 
corresponds to a convergent Jacobi-Perron fraction,  and
let $W=F(V)$.    Our main result is as follows. 
\begin{theorem}\label{thm1}
The set $V$ is a generic subset of $T(g)$ and  map $F$
has the following properties:
(i) $V\cong W\times (0,\infty)$ is a trivial fiber bundle, whose projection map
$\pi: V\to W$ coincides with $F$ and 
(ii) $F$ is a covariant functor which maps isomorphic Riemann surfaces  $S,S'\in V$ 
to the stably isomorphic toric $AF$-algebras  ${\Bbb A}_{\theta}, {\Bbb A}_{\theta'} \in W$.
\end{theorem}
The article is organized as follows.  Preliminary facts are reviewed in 
Section 2.   Theorem  \ref{thm1}  is proved in Section 3.

\section{Preliminaries}
\subsection{Measured foliations and $T(g)$}
 A measured foliation, ${\Bbb F}$, on a surface $X$
is a  partition of $X$ into the singular points $x_1,\dots,x_n$ of
order $k_1,\dots, k_n$ and regular leaves (1-dimensional submanifolds). 
On each  open cover $U_i$ of $X-\{x_1,\dots,x_n\}$ there exists a non-vanishing
real-valued closed 1-form $\phi_i$  such that: 
(i)  $\phi_i=\pm \phi_j$ on $U_i\cap U_j$;
(ii) at each $x_i$ there exists a local chart $(u,v):V\to {\Bbb R}^2$
such that for $z=u+iv$, it holds $\phi_i=Im~(z^{k_i\over 2}dz)$ on
$V\cap U_i$ for some branch of $z^{k_i\over 2}$. 
The pair $(U_i,\phi_i)$ is called an atlas for measured foliation ${\Bbb F}$.
Finally, a measure $\mu$ is assigned to each segment $(t_0,t)\in U_i$, which is  transverse to
the leaves of ${\Bbb F}$, via the integral $\mu(t_0,t)=\int_{t_0}^t\phi_i$. The 
measure is invariant along the leaves of ${\Bbb F}$, hence the name. 
We refer the reader to [Thurston 1988] \cite{Thu1}  and [Fathi, Laudenbach \& Po\'enaru 1979]
\cite{FLP} for a systematic account  of measured foliations.

Let $S$ be a Riemann surface, and $q\in H^0(S,\Omega^{\otimes 2})$ a holomorphic
quadratic differential on $S$. The lines $Re~q=0$ and $Im~q=0$ define a pair
of measured foliations on $R$, which are transversal to each other outside the set of 
singular points. The set of singular points is common to both foliations and coincides
with the zeroes of $q$. The above measured foliations are said to represent the  vertical and horizontal 
 trajectory structure  of $q$, respectively. 
Let $T(g)$ be the Teichm\"uller space of the topological surface $X$ of genus $g\ge 1$,
i.e. the space of the complex structures on $X$. 
Consider the vector bundle $p: Q\to T(g)$ over $T(g)$ whose fiber above a point 
$S\in T(g)$ is the vector space $H^0(S,\Omega^{\otimes 2})$.   
Given non-zero $q\in Q$ above $S$, we can consider horizontal measured foliation
${\Bbb F}_q\in \Phi_X$ of $q$, where $\Phi_X$ denotes the space of equivalence
classes of measured foliations on $X$. If $\{0\}$ is the zero section of $Q$,
the above construction defines a map $Q-\{0\}\longrightarrow \Phi_X$. 
For any ${\Bbb F}\in\Phi_X$, let $E_{\Bbb F}\subset Q-\{0\}$ be the fiber
above ${\Bbb F}$. In other words, $E_{\Bbb F}$ is a subspace of the holomorphic 
quadratic forms whose horizontal trajectory structure coincides with the 
measured foliation ${\Bbb F}$. 
Note that, if ${\Bbb F}$ is a measured foliation with the simple zeroes (a generic case),  
then $E_{\Bbb F}\cong {\Bbb R}^n - 0$, while $T(g)\cong {\Bbb R}^n$, where $n=6g-6$ if
$g\ge 2$ and $n=2$ if $g=1$.
\begin{lemma}
{\bf  [Hubbard \& Masur  1979]  \cite{HuMa1}}             
The restriction of $p$ to $E_{\Bbb F}$ defines a homeomorphism (an embedding)
$h_{\Bbb F}: E_{\Bbb F}\to T(g)$.
\end{lemma}
The Hubbard-Masur result implies that the measured foliations  parametrize  
the space $T(g)-\{pt\}$, where $pt= h_{\Bbb F}(0)$.
Indeed, denote by  ${\Bbb F}'$ a vertical trajectory structure of  $q$. Since ${\Bbb F}$
and ${\Bbb F}'$ define $q$, and ${\Bbb F}=Const$ for all $q\in E_{\Bbb F}$, one gets a homeomorphism 
between $T(g)-\{pt\}$ and  $\Phi_X$, where $\Phi_X\cong {\Bbb R}^n-0$ is the space of 
equivalence classes of the measured foliations ${\Bbb F}'$ on $X$. 
Note that the above parametrization depends on a foliation ${\Bbb F}$.
However, there exists a unique canonical homeomorphism $h=h_{\Bbb F}$
as follows. Let $Sp ~(S)$ be the length spectrum of the Riemann surface
$S$ and $Sp ~({\Bbb F}')$ be the set positive reals $\inf \mu(\gamma_i)$,
where $\gamma_i$ runs over all simple closed curves, which are transverse 
to the foliation ${\Bbb F}'$. A canonical homeomorphism 
$h=h_{\Bbb F}: \Phi_X\to T(g)-\{pt\}$ is defined by the formula
$Sp ~({\Bbb F}')= Sp ~(h_{\Bbb F}({\Bbb F}'))$ for $\forall {\Bbb F}'\in\Phi_X$. 
Thus, the following corollary is true.
\begin{corollary}\label{cor1}
There exists a unique homeomorphism $h:\Phi_X\to T(g)-\{pt\}$.
\end{corollary}
Recall that $\Phi_X$ is the space of equivalence classes of measured
foliations on the topological surface $X$. Following [Douady \& Hubbard 1975]
\cite{DoHu1}, we consider a coordinate system on $\Phi_X$,
suitable for the proof of theorem \ref{thm1}. 
For clarity, let us make a generic assumption that $q\in H^0(S,\Omega^{\otimes 2})$
is a non-trivial holomorphic quadratic differential with only simple zeroes. 
We wish to construct a Riemann surface of $\sqrt{q}$, which is a double cover
of $S$ with ramification over the zeroes of $q$. Such a surface, denoted by
$\widetilde S$, is unique and has an advantage of carrying a holomorphic
differential $\omega$, such that $\omega^2=q$. We further denote by 
$\pi:\widetilde S\to S$ the covering projection. The vector space
$H^0(\widetilde S,\Omega)$ splits into the direct sum
$H^0_{even}(\widetilde S,\Omega)\oplus H^0_{odd}(\widetilde S,\Omega)$
in view of  the involution $\pi^{-1}$ of $\widetilde S$, and
the vector space $H^0(S,\Omega^{\otimes 2})\cong H^0_{odd}(\widetilde S,\Omega)$.
Let $H_1^{odd}(\widetilde S)$ be the odd part of the homology of $\widetilde S$
relatively  the zeroes of $q$.   Consider the pairing
$H_1^{odd}(\widetilde S)\times H^0(S, \Omega^{\otimes 2})\to {\Bbb C}$,
defined by the integration  $(\gamma, q)\mapsto \int_{\gamma}\omega$. 
We shall take the associated map
$\psi_q: H^0(S,\Omega^{\otimes 2})\to Hom~(H_1^{odd}(\widetilde S); {\Bbb C})$
and let $h_q= Re~\psi_q$. 
\begin{lemma}{\bf [Douady \& Hubbard 1975]  \cite{DoHu1}}
The map
\begin{equation}\label{eq3a}
h_q: H^0(S, \Omega^{\otimes 2})\longrightarrow Hom~(H_1^{odd}(\widetilde S); {\Bbb R})
\end{equation}
is an ${\Bbb R}$-isomorphism. 
\end{lemma}
Since  each  ${\Bbb F}\in \Phi_X$ is the  vertical foliation 
$Re~q=0$ for a $q\in H^0(S, \Omega^{\otimes 2})$, the Douady-Hubbard lemma
implies that $\Phi_X\cong Hom~(H_1^{odd}(\widetilde S); {\Bbb R})$.
By  formulas for the relative homology, 
one finds that $H_1^{odd}(\widetilde S)\cong {\Bbb Z}^{n}$,
where $n=6g-6$ if $g\ge 2$ and $n=2$ if $g=1$.     
Finally, each $h\in Hom~({\Bbb Z}^{n}; {\Bbb R})$ is given
by the reals  $\lambda_1=h(e_1),\dots, \lambda_{n}=h(e_{n})$,
where $(e_1,\dots, e_{n})$ is a basis in ${\Bbb Z}^{n}$.
The numbers   $(\lambda_1,\dots,\lambda_{n})$ are the coordinates in the space $\Phi_X$
and, in view of the corollary \ref{cor1}, in  the Teichm\"uller space $T(g)$.

\subsection{The Jacobi-Perron continued fraction}
Let $a_1,a_2\in {\Bbb N}$ such that $a_2\le a_1$. Recall that the greatest common
divisor of $a_1,a_2$, $GCD(a_1,a_2)$, can be determined from the Euclidean algorithm:
\begin{equation}\label{eq5}
\left\{
\begin{array}{cc}
a_1 &= a_2b_1 +r_3\nonumber\\
a_2 &= r_3b_2 +r_4\nonumber\\
r_3 &= r_4b_3 +r_5\nonumber\\
\vdots & \nonumber\\
r_{k-3} &= r_{k-2}b_{k-1}+r_{k-1}\nonumber\\
r_{k-2} &= r_{k-1}b_k,
\end{array}
\right.
\end{equation}
where $b_i\in {\Bbb N}$ and $GCD(a_1,a_2)=r_{k-1}$. 
The Euclidean algorithm can be written as the regular continued 
fraction
\begin{equation}\label{eq6a}
\theta={a_1\over a_2}=b_1+{1\over\displaystyle b_2+
{\strut 1\over\displaystyle +\dots+ {1\over b_k}}}
=[b_1,\dots b_k].
\end{equation}
If $a_1, a_2$ are non-commensurable, in the sense that $\theta\in {\Bbb R}-{\Bbb Q}$,
then the Euclidean algorithm never stops and $\theta=[b_1, b_2, \dots]$. Note that the regular  
continued fraction can be written in the matrix form:
\begin{equation}\label{eq7}
\left(
\begin{matrix}
1\cr \theta
\end{matrix}
\right)=
\lim_{k\to\infty} \left(
\begin{matrix}
0 & 1\cr 1 & b_1
\end{matrix}
\right)\dots
\left(
\begin{matrix}
0 & 1\cr 1 & b_k
\end{matrix}
\right)
\left(
\begin{matrix}
0\cr 1
\end{matrix}
\right). 
\end{equation}
The Jacobi-Perron algorithm and connected (multidimensional) continued 
fraction generalizes the Euclidean algorithm to the case $GCD(a_1,\dots,a_n)$
when $n\ge 2$. Namely, let $\lambda=(\lambda_1,\dots,\lambda_n)$,
$\lambda_i\in {\Bbb R}-{\Bbb Q}$ and  $\theta_{i-1}={\lambda_i\over\lambda_1}$, where
$1\le i\le n$.   The continued fraction 
$$
\left(
\begin{matrix}
1\cr \theta_1\cr\vdots\cr\theta_{n-1}
\end{matrix} 
\right)=
\lim_{k\to\infty} 
\left(
\begin{matrix}
0 &  0 & \dots & 0 & 1\cr
              1 &  0 & \dots & 0 & b_1^{(1)}\cr
              \vdots &\vdots & &\vdots &\vdots\cr
              0 &  0 & \dots & 1 & b_{n-1}^{(1)}
\end{matrix}              
\right)
\dots 
\left(
\begin{matrix}
0 &  0 & \dots & 0 & 1\cr
              1 &  0 & \dots & 0 & b_1^{(k)}\cr
              \vdots &\vdots & &\vdots &\vdots\cr
              0 &  0 & \dots & 1 & b_{n-1}^{(k)}
              \end{matrix}
              \right)
\left(
\begin{matrix}
0\cr 0\cr\vdots\cr 1
\end{matrix}
\right),
$$
where $b_i^{(j)}\in {\Bbb N}\cup\{0\}$, is called the {\it Jacobi-Perron
algorithm (JPA)}. Unlike the regular continued fraction algorithm,
the JPA may diverge for certain vectors $\lambda\in {\Bbb R}^n$. However, 
for points of a generic subset of ${\Bbb R}^n$, the JPA converges. 
 The convergence of the JPA algorithm can be characterized in terms of
the measured foliations. Let ${\Bbb F}\in\Phi_X$ be a measured foliation
on the surface $X$ of genus $g\ge 1$. Recall that ${\Bbb F}$ is called uniquely
ergodic if every invariant measure of ${\Bbb F}$ is a multiple
of the Lebesgue measure. It is known that there exists
a generic subset $V\subset \Phi_X$ such that each ${\Bbb F}\in V$
is uniquely ergodic [Masur 1982]  \cite{Mas1} and [Veech 1982]   \cite{Vee1}.
We let $\lambda=(\lambda_1,\dots,\lambda_{n})$ be the vector with
coordinates $\lambda_i=\mu ({\gamma_i})$, where $\gamma_i\in H_1^{odd}(\widetilde S)$;
by an abuse of notation,  we shall say that $\lambda\in V$. 
In view of a bijection  between measured foliations and the interval exchange transformations ([Masur 1982]  \cite{Mas1}), 
the following characterization of convergence of the JPA is true.
\begin{lemma}\label{lm0}
{\bf [Bauer 1996]  \cite{Bau1}}
The JPA converges if and only if $\lambda\in V\subset {\Bbb R}^{n}$.  
\end{lemma}

\section{Proof of theorem 1}
Let us outline the proof. We shall  consider the following sets of objects:

\medskip
(i) generic Riemann surfaces $V$;

\smallskip
(ii) pseudo-lattices   ${\Bbb PL}$,   see  [Manin 2004 ]  \cite{Man2};

\smallskip
(iii) projective pseudo-lattices ${\Bbb PPL}$;

\smallskip
 (iv)  toric $AF$-algebras $W$.

\bigskip\noindent
The proof takes  the following steps:

\medskip
(a) to show that  $V\cong {\Bbb PL}$
are equivalent categories such that isomorphic Riemann surfaces 
$S,S'\in V$ map to isomorphic pseudo-lattices
$PL, PL'\in {\Bbb PL}$;

\smallskip
(b) a non-injective functor $F: {\Bbb PL}\to {\Bbb PPL}$
is constructed. The $F$ maps isomorphic pseudo-lattices to isomorphic
projective pseudo-lattices and $Ker ~F\cong (0,\infty)$;

\smallskip
(c) to show that a subcategory $U\subseteq {\Bbb PPL}$ and  $W$ are the 
equivalent categories.

\vskip0.5cm\noindent
In other words, we have the following diagram:
\begin{equation}\label{eq1}
V
\buildrel\rm\alpha
\over\longrightarrow
{\Bbb PL}
\buildrel\rm F
\over\longrightarrow
U
\buildrel\rm \beta
\over\longrightarrow
 W,
\end{equation}
where $\alpha$ is an injective map, $\beta$ is a bijection  and $Ker ~F\cong (0,\infty)$.

\bigskip\noindent
{\bf (i) Category $V$}. 
A {\it Riemann surface}  is a triple $(X, S, j)$, where
$X$ is a topological surface of genus $g\ge 1$,  $j: X\to S$ is a
complex (conformal) parametrization of $X$ and $S$ is a Riemann surface.
A {\it morphism} of Riemann surfaces  $(X, S, j)\to (X, S', j')$   
is a biholomorphic map modulo the ones, which are isotopic to the identity map
with respect to a fixed topological marking of $X$.  
A  category of generic Riemann surfaces $V$  consists of $Ob~({\Bbb S})$
which are Riemann surfaces  $S\in V\subset T(g)$ and morphisms $H(S,S')$
between $S,S'\in Ob~(V)$ which coincide with the morphisms
specified above. For any  $S,S',S''\in Ob~({\Bbb S})$
and any morphisms $\varphi': S\to S'$, $\varphi'': S'\to S''$ a 
morphism $\phi: S\to S''$ is the composite of $\varphi'$ and
$\varphi''$, which we write as $\phi=\varphi''\varphi'$. 
The identity  morphism, $1_S$, is a morphism $H(S,S)$.

\bigskip\noindent
{\bf (ii)  Category ${\Bbb PL}$}. 
A {\it pseudo-lattice} (of rank $n$) is a triple $(\Lambda, {\Bbb R}, j)$, where
$\Lambda\cong {\Bbb Z}^n$ and $j: \Lambda\to {\Bbb R}$ is a homomorphism.  
A morphism of pseudo-lattices $(\Lambda, {\Bbb R}, j)\to (\Lambda, {\Bbb R}, j')$
is a commutative diagram:

\begin{picture}(300,110)(-80,-5)
\put(40,70){\vector(0,-1){35}}
\put(120,70){\vector(0,-1){35}}
\put(65,23){\vector(1,0){33}}
\put(65,83){\vector(1,0){33}}
\put(38,20){${\Bbb Z}^n$}
\put(118,20){$\Bbb R$}
\put(37,80){${\Bbb Z}^n$}
\put(117,80){${\Bbb R}$}
\put(50,50){$\varphi$}
\put(130,50){$\psi$}
\put(77,28){$j'$}
\put(77,88){$j$}
\end{picture}

\noindent
where $\varphi$ is a group homomorphism and $\psi$ is an inclusion 
map, i.e. $j'(\Lambda')\subseteq j(\Lambda)$.  
Any isomorphism class of a pseudo-lattice contains
a representative given by $j: {\Bbb Z}^n\to  {\Bbb R}$ such that 
$$j(1,0,\dots, 0)=\lambda_1, \quad j(0,1,\dots,0)=\lambda_2,\quad \dots, \quad  j(0,0,\dots,1)=\lambda_n,$$
where $\lambda_1,\lambda_2,\dots,\lambda_n$ are positive reals.
The pseudo-lattices of rank $n$  make up a category, which we denote by ${\Bbb PL}_n$.

\bigskip
The following lemma  says that the  ${\Bbb Z}$-module ${\Bbb Z}\lambda_1+\dots+
{\Bbb Z}\lambda_n$ is an  invariant of the isomorphism class of the Riemann surface $S$; 
in other words,  the action of  mapping class group $Mod~(X)$ on such a module corresponds to a transformation of
 basis of the module. 
\begin{lemma}\label{lm1}
Let $g\ge 2$ (resp. $g=1$) and $n=6g-6$ (resp. $n=2$). 
There exists an injective covariant functor $\alpha:  V\to {\Bbb PL}_n$
which maps isomorphic Riemann surfaces  $S,S'\in V$ to the isomorphic
pseudo-lattices $PL,PL'\in {\Bbb PL}_n$.
\end{lemma}
{\it Proof.} Let $\alpha: T(g)-\{pt\}\to Hom~(H_1^{odd}(\widetilde S); {\Bbb R})-0$
be a Hubbard-Masur map. Since $\alpha$ is a homeomorphism 
between the respective spaces, we conclude that $\alpha$ is an injective
map. The first claim of lemma is proved.

Let us show that $\alpha$ sends morphisms of ${\Bbb S}$ to morphisms
of ${\Bbb PL}$. Let $\varphi\in Mod~(X)$ be a diffeomorphism of $X$.
Suppose that all the zeroes of measured foliations are generic (simple)
and  let $p:\widetilde X\to X$ be the double cover of $X$. 
(Note that the case of torus does not require
a double cover, and thus one can assert $p=Id$ in the argument below.)
 Denote by $\widetilde\varphi$ a diffeomorphism
of $\widetilde X$, which makes the following diagram commutative:

\begin{picture}(300,110)(-80,-5)
\put(40,70){\vector(0,-1){35}}
\put(120,70){\vector(0,-1){35}}
\put(65,23){\vector(1,0){33}}
\put(65,83){\vector(1,0){33}}
\put(38,20){$X$}
\put(118,20){$X$}
\put(37,80){$\widetilde X$}
\put(117,80){$\widetilde X$}
\put(50,50){$p$}
\put(130,50){$p$}
\put(77,28){$\varphi$}
\put(77,88){$\widetilde\varphi$}
\end{picture}

\noindent
One can consider the effect of $\varphi, \widetilde\varphi$ and $p$
on the respective (relative) integral homology groups:

\begin{picture}(300,110)(-40,-5)
\put(40,70){\vector(0,-1){35}}
\put(180,70){\vector(0,-1){35}}
\put(95,23){\vector(1,0){33}}
\put(95,83){\vector(1,0){33}}
\put(10,20){$H_1(X, Sing~{\Bbb F})$}
\put(148,20){$H_1(X, Sing~{\Bbb F})$}
\put(-20,80){$H_1^{odd}(\widetilde X)\oplus H_1^{even}(\widetilde X)$}
\put(147,80){$H_1^{odd}(\widetilde X)\oplus H_1^{even}(\widetilde X)$}
\put(50,50){$p_*$}
\put(190,50){$p_*$}
\put(107,28){$\varphi_*$}
\put(107,88){$\widetilde\varphi_*$}
\end{picture}

\noindent
where $Ker~p_*\cong H_1^{even}(\widetilde X)$.  
Since $p_*: H_1^{odd}(\widetilde X)\to H_1(X, Sing~{\Bbb F})$
is an isomorphism, we conclude that $\widetilde\varphi_*\in GL_n({\Bbb Z})$,
where $n=dim~H_1^{odd}(\widetilde X)$. It is easy to see, that  
 $\widetilde\varphi_*$ acts on a pseudo-lattice by a transformation
of its  basis, and  therefore, $\widetilde\varphi_*\in Mor~({\Bbb PL})$.

\medskip
Let us show that $\alpha$ is a functor.  Indeed, let $S,S'\in V$ be isomorphic
Riemann surfaces, such that $S'=\varphi(S)$ for a $\varphi\in Mod~(X)$. 
Let $a_{ij}$ be the elements of matrix $\widetilde\varphi_*\in GL_n({\Bbb Z})$. 
Recall that: 
\begin{equation}\label{eq2}
\lambda_i=\int_{\gamma_i}\phi
\end{equation}
for a closed 1-form $\phi= Re~\omega$
and $\gamma_i\in H_1^{odd}(\widetilde X)$. Then
\begin{equation}\label{eq3}
\gamma_j=\sum_{i=1}^n a_{ij}\gamma_i, \qquad j=1,\dots, n, 
\end{equation}
are the elements of a new basis in $H_1^{odd}(\widetilde X)$. By the integration rules,
\begin{equation}\label{eq4}
\lambda_j'= \int_{\gamma_j}\phi=
\int_{\sum a_{ij}\gamma_i}\phi=
\sum_{i=1}^n a_{ij}\lambda_i.
\end{equation}
Finally, let $j(\Lambda)={\Bbb Z}\lambda_1+\dots+{\Bbb Z}\lambda_n$
and $j'(\Lambda)={\Bbb Z}\lambda_1'+\dots+{\Bbb Z}\lambda_n'$.
Since $\lambda_j'= \sum_{i=1}^n a_{ij}\lambda_i$ and $(a_{ij})\in GL_n({\Bbb Z})$,
we conclude that: 
\begin{equation}\label{eq5b}
   j(\Lambda)=j'(\Lambda).
\end{equation}
In other words, the ${\Bbb Z}$-module ${\Bbb Z}\lambda_1+\dots+{\Bbb Z}\lambda_n$
is an invariant of $Mod~(X)$.  In particular, the pseudo-lattices
$(\Lambda, {\Bbb R}, j)$ and $(\Lambda, {\Bbb R}, j')$ are isomorphic. Hence, 
$\alpha: V\to {\Bbb PL}$  maps isomorphic Riemann surfaces to
the isomorphic pseudo-lattices, i.e. $\alpha$ is a functor.

\medskip
Finally, let us show that $\alpha$ is a covariant functor. 
Indeed, let $\varphi_1,\varphi_2\in Mor ({\Bbb S})$. Then  
$\alpha(\varphi_1\varphi_2)= (\widetilde{\varphi_1\varphi_2})_*=(\widetilde\varphi_1)_*(\widetilde\varphi_2)_*=
\alpha(\varphi_1)\alpha(\varphi_2)$. Lemma \ref{lm1} follows.
$\square$

\bigskip\noindent
{\bf (iii)  Category ${\Bbb PPL}$}. 
A  {\it projective pseudo-lattice} (of rank $n$) is a triple 
$(\Lambda, {\Bbb R}, j)$, where $\Lambda\cong {\Bbb Z}^n$ and $j: \Lambda\to {\Bbb R}$ is a
homomorphism. A morphism of projective pseudo-lattices
$(\Lambda, {\Bbb C}, j)\to (\Lambda, {\Bbb R}, j')$
is a commutative diagram:

\begin{picture}(300,110)(-80,-5)
\put(40,70){\vector(0,-1){35}}
\put(120,70){\vector(0,-1){35}}
\put(65,23){\vector(1,0){33}}
\put(65,83){\vector(1,0){33}}
\put(38,20){${\Bbb Z}^n$}
\put(118,20){$\Bbb R$}
\put(37,80){${\Bbb Z}^n$}
\put(117,80){${\Bbb R}$}
\put(50,50){$\varphi$}
\put(130,50){$\psi$}
\put(77,28){$j'$}
\put(77,88){$j$}
\end{picture}

\noindent
where $\varphi$ is a group homomorphism and $\psi$ is an  ${\Bbb R}$-linear
map. (Notice that unlike the case of pseudo-lattices, $\psi$ is a scaling
map as opposite to an inclusion map. This allows to the two pseudo-lattices
to be projectively equivalent, while being distinct in the category ${\Bbb PL}_n$.) 
It is not hard to see that any isomorphism class of a projective pseudo-lattice 
contains a representative given by $j: {\Bbb Z}^n\to  {\Bbb R}$ such that
$$j(1,0,\dots,0)=1,\quad
j(0,1,\dots,0)=\theta_1,\quad  \dots, \quad  j(0,0,\dots,1)=\theta_{n-1},$$
  where $\theta_i$ are  positive reals.
The projective pseudo-lattices of rank $n$  make up a category, which we denote by ${\Bbb PPL}_n$.

\bigskip\noindent
{\bf (iv) Category $W$}. 
Let $\theta=(\theta_1,\dots,\theta_{n-1})$.  Then
toric $AF$-algebras ${\Bbb A}_{\theta}$ make a category; 
morphisms in the category are  stable isomorphism between 
toric $AF$-algebras.    We shall denote such a category by $W_n$.  
\begin{lemma}\label{lm2}
Let $U_n\subseteq {\Bbb PPL}_n$ be a subcategory consisting of the projective
pseudo-lattices $PPL=PPL(1,\theta_1,\dots,\theta_{n-1})$ for which the Jacobi-Perron
fraction of the vector $(1,\theta_1,\dots,\theta_{n-1})$ converges to the vector.
Define a map  $\beta: U_n\to W_n$  by the formula
$PPL(1,\theta_1,\dots,\theta_{n-1})\mapsto {\Bbb A}_{\theta}.$
Then $\beta$ is a bijective functor,
which maps isomorphic projective pseudo-lattices to the stably isomorphic 
toric $AF$-algebras. 
 \end{lemma}
{\it Proof.} It is evident that $\beta$ is injective and surjective. Let
us show that $\beta$ is a functor. Indeed,   every totally 
ordered abelian group of rank $n$ has form ${\Bbb Z}+\theta_1 {\Bbb Z}+\dots+ {\Bbb Z}\theta_{n-1}$,
see e.g. [Effros 1981]   \cite{E}, Corollary 4.7.
The latter   is a projective pseudo-lattice $PPL$  from the category $U_n$. 
On the other hand,  each  $PPL$ defines a stable isomorphism class of 
the $AF$-algebra ${\Bbb A}_{\theta_1,\dots,\theta_{n-1}}\in W_n$  [Elliott 1976] \cite{Ell1}. 
Therefore, $\beta$  maps isomorphic projective pseudo-lattices (from the set $U_n$) to the stably isomorphic toric 
$AF$-algebras,  and {\it vice versa}. Lemma \ref{lm2} follows.  
 $\square$

\bigskip
Let $PL(\lambda_1,\lambda_2,\dots, \lambda_n)\in {\Bbb PL}_n$ and 
$PPL(1,\theta_1,\dots,\theta_{n-1})\in {\Bbb PPL}_n$. 
To finish the proof of theorem \ref{thm1}, it remains to show the following.
\begin{lemma}\label{lm3}
Let $F: {\Bbb PL}_n\to {\Bbb PPL}_n$ be a map given by formula
$$PL(\lambda_1,\lambda_2,\dots, \lambda_n)\mapsto PPL\left(1, {\lambda_2\over\lambda_1},\dots, {\lambda_n\over\lambda_1}\right).$$
Then $Ker~F=(0,\infty)$ and  $F$ is a functor which maps  isomorphic pseudo-lattices to isomorphic projective 
pseudo-lattices. 
\end{lemma}
{\it Proof.} Indeed, $F$ can be thought as   a map from ${\Bbb R}^n$ to ${\Bbb R}P^{n-1}$. Hence 
$Ker~F= \{\lambda_1 : \lambda_1>0\}\cong (0,\infty)$. The second part of lemma is evident.  
$\square$

\medskip
Assuming $n=6g-6$ (resp. $n=2$) for $g\ge 2$ (resp. $g=1$), one gets items (i) and (ii) of the second part of 
theorem \ref{thm1} from lemmas \ref{lm1}-\ref{lm3};   the first part of theorem \ref{thm1} (i.e. that $V$
is generic) follows from lemma \ref{lm0}. 
$\square$

\bibliographystyle{amsplain}

\begin{thebibliography}{99}



\bibitem{Bau1}
M.~Bauer, \textit{A characterization of uniquely ergodic interval exchange
maps in terms of the Jacobi-Perron algorithm}, 
Bol. Soc. Bras. Mat. {\bf 27} (1996),
109-128. 


\bibitem{B}
L.~Bernstein, \textit{The Jacobi-Perron Algorithm, its Theory and Applications},
Lect. Notes in Math. {\bf 207}, Springer 1971. 


\bibitem{Bra1}
O.~Bratteli, \textit{Inductive limits of finite dimensional
$C^*$-algebras}, Trans. Amer. Math. Soc. {\bf 171} (1972),
195-234. 




\bibitem{DoHu1}
A.~Douady and J.~Hubbard, \textit{On the density of Strebel differentials},
Inventiones Math. {\bf 30} (1975), 175-179. 


\bibitem{E}
E.~G.~Effros, \textit{Dimensions and $C^*$-Algebras}, in: Conf. Board of the Math.
Sciences, Regional conference series in Math., No. {\bf 46}, AMS,  1981.



\bibitem{EfSh2}
E.~G.~Effros and C.~L.~Shen,
\textit{Approximately finite $C^*$-algebras and continued fractions},
Indiana Univ. Math. J. {\bf 29} (1980), 191-204.



\bibitem{Ell1}
G.~A.~Elliott, \textit{On the classification of inductive limits of sequences
of semisimple finite-dimensional algebras}, J.~Algebra {\bf 38} (1976), 29-44. 


\bibitem{FLP}
A.~Fathi, F.~Laudenbach and V.~Po\'enaru,
\textit{Travaux de Thurston sur les Surfaces},  Ast\'erisque 
{\bf 66-67},  1979.



\bibitem{HuMa1} 
J.~Hubbard and H.~Masur, \textit{Quadratic differentials and foliations}, 
Acta Math. {\bf 142} (1979), 221-274.


\bibitem{Kle1}
F.~Klein, \textit{Ausgew\"alte Kapitel der Zahlentheorie}, 
Math. Annalen {\bf 48} (1896/97). 












\bibitem{Man2}
Yu.~I.~Manin, \textit{Real multiplication and noncommutative geometry},
in ``Legacy of Niels Hendrik Abel'', 685-727, Springer, 2004.



\bibitem{Mas1}
H.~Masur, \textit{Interval exchange transformations and measured foliations},
Ann. of Math. (2) {\bf 115} (1982), no. 1, 169--200.












\bibitem{S}
K.~Strebel, \textit{Quadratic Differentials}, Springer, 1984. 



\bibitem{Thu1}
W.~P.~Thurston, \textit{On the geometry and dynamics of diffeomorphisms
of surfaces},  Bull. Amer. Math. Soc. {\bf 19} (1988), 417-431.


    
        
\bibitem{Vee1}
W.~A.~Veech, \textit{Gauss measures for transformations on the space of 
interval exchange maps},  Ann. of Math. (2) {\bf 115}  (1982), no. 1, 201-242.
    


\end{thebibliography}

\end{document}